\documentclass{daj}

\usepackage[english]{babel}
\usepackage[applemac]{inputenc}
\usepackage[T1]{fontenc}
\usepackage{palatino}
\usepackage{bbm}
\usepackage{cite}
\usepackage{float}
\usepackage{amsmath}
\usepackage{amssymb}
\usepackage{amsthm}
\usepackage{amsfonts}
\usepackage{graphicx}
\usepackage{mathtools}
\usepackage{enumitem}
\usepackage{color}

\usepackage{xcolor}
\usepackage[pagebackref]{hyperref}
\hypersetup{
   colorlinks,
}
\usepackage[margin=3cm]{geometry}

\mathtoolsset{showonlyrefs}

\newcommand{\R}{\mathbb{R}}

\newcommand{\eps}{\varepsilon}
\newcommand{\DK}{\mathcal{DK}}

\DeclareMathOperator{\good}{\textup{\textbf{Good}}}
\DeclareMathOperator{\bad}{\textup{\textbf{Bad}}}

\newcommand{\be}{\begin{equation}}
\newcommand{\ee}{\end{equation}}

\DeclareMathOperator{\Haus}{\dim_H}

\DeclareMathOperator{\Assouad}{\dim_A}

\DeclareMathOperator{\uboxd}{\overline{\dim}_B}
\DeclareMathOperator{\lboxd}{\underline{\dim}_B}
\DeclareMathOperator{\dir}{dir}
\DeclareMathOperator{\dist}{dist}

\newtheorem{thm}{Theorem}[section]

\newtheorem{prob}[thm]{Problem}

\newtheorem*{thm*}{Theorem}
\newtheorem*{conj*}{Conjecture}
\newtheorem*{lma*}{Lemma}

\dajAUTHORdetails{%
  title = {On sets containing a unit distance in every direction}, 
  author = {Pablo Shmerkin and Han Yu},
  plaintextauthor = {Pablo Shmerkin, Han yu},
  keywords = {Box dimension, unit distances, difference sets},
}   

\dajEDITORdetails{%
   year={2021},
   number={5},
   received={10 August 2020},   
   revised={14 March 2021},    
   published={5 July 2021},  
   doi={10.19086/da.22058},       
}   

\begin{document}

\begin{frontmatter}[classification=text]


\author[pablo]{Pablo Shmerkin \thanks{PS was supported by Royal Society grant IES\textbackslash R1\textbackslash 191195 and by a St Andrews Global Fellowship}}
\author[han]{Han Yu\thanks{HY has received funding from the European Research Council (ERC) under the European Union's Horizon 2020 research and innovation programme (grant agreement No. 803711). }}

\begin{abstract}
We investigate the box dimensions of compact sets in $\mathbb{R}^2$ that contain a unit distance in every direction (such sets may have zero Hausdorff dimension). Among other results, we show that the lower box dimension must be at least $\frac{4}{7}$ and can be as low as $\frac{2}{3}$. This quantifies in a certain sense how far the unit circle is from being a difference set.
\end{abstract}
\end{frontmatter}


\section{Introduction}
In this paper, we are interested in the size of sets $A\subset\R^2$ such that $A-A\supset S^1$ or, in other words, such that for all $e\in S^1$ there are $x,y\in A$ such that $y-x=e$. Thus, such sets can be seen as a variant of Kakeya sets which contain, instead of a unit segment in every direction, just the endpoints of the segment. They can also be seen as sets that contain many unit distances.  We call bounded sets $A\subset\R^2$ such that $A-A\supset S^{1}$ \emph{dipole Kakeya sets} and denote the collection of all of them by $\DK_2$. Similarly, we denote by $\DK_n$ the collection of bounded sets $A\subset\R^n$ such that $A-A\supset S^{n-1}$.

A first attempt to quantify the size of dipole Kakeya sets might be via Hausdorff dimension. However, using the Baire category method in \cite{CCHK18} one can show that there are compact Kakeya dipole sets of zero Hausdorff dimension; in fact they are `typically' of zero Hausdorff dimension in the appropriate context. Therefore, in this paper we focus on the box dimensions of dipole Kakeya sets (see Section \ref{Pre} for definitions of dimensions); this problem was suggested to us by Alan Chang. Recall that
\[
\lboxd(A-A) \le \lboxd(A\times A) = 2\lboxd(A),
\]
where the left-hand inequality size follows since $(x,y)\mapsto y-x$ is Lipschitz and the right-hand side equality is immediate from the definition. Hence, dipole Kakeya sets in $\mathbb{R}^n$ have (lower) box dimension at least $(n-1)/2$, and any improvement on this bound quantifies, in some sense, the difficulty that the unit sphere has in behaving like a difference set. The following is our main result:
\begin{thm}\label{MAIN}
	 Let $A\subset\mathbb{R}^2$ be a (bounded) dipole Kakeya set. Then we have
\[
	\lboxd A\geq \frac{4}{7},
	\]
	and
\[
    \Assouad A \geq \frac{2}{3}.
\]
where $\Assouad$ denotes Assouad dimension.

On the other hand, for each $n\ge 2$ there exist compact dipole Kakeya sets in $\mathcal{DK}_n$ with $\lboxd(A)\le n(n-1)/(2n-1)$.
\end{thm}
As a remark, we point out that the examples we construct in this paper also have zero Hausdorff dimension, see \S\ref{subsec:Hausdorff}.

The longstanding Erd\H{o}s unit distance conjecture says that for a finite planar set $A$ with $k\geq 2$ points, the number of pairs of points in $A$ at distance $1$ from each other is bounded above by $C_\eps k^{1+\eps}$ for each $\eps>0$. The best known bound, from \cite{SST84},  is $C_{1/3}k^{4/3}$ (that is, $\eps=1/3$). By analogy with this discrete problem, one might expect that a dipole Kakeya set has dimension at least $3/4.$ As we have seen, this is not true for lower box dimension, but it might be true for upper box dimension. In support of this, in Section \ref{RE} we construct a set in $\mathcal{DK}_2$ with upper box dimension at most $3/4$. We pose here the following problem in $\mathbb{R}^2$.
\begin{prob}\label{Con}
Is it true that
\[
\inf\{ \uboxd(A): A\subset\R^2, A-A \supset S^1 \} = \frac{3}{4}
\]
and
\[
\inf\{ \lboxd(A): A\subset\R^2, A-A\supset S^1 \} = \frac{2}{3} \, \text{?}
\]
\end{prob}

\section{Preliminaries}\label{Pre}
In this section, we shall give a brief introduction to the notions of dimensions used in this paper. In particular, we shall encounter Hausdorff dimension, box dimensions and the Assouad dimension. For more detailed background, see \cite[Chapters 2,3]{Falconer14}, \cite[Chapters 4,5]{Mattila95} and \cite{Fraser14},

\subsection{Hausdorff dimension}

Given an exponent $s\in\mathbb{R}^+$, define the $s$-Hausdorff content of a set $F\subset\R^n$ as
\[
\mathcal{H}^s_\infty(F)=\inf\left\{\sum_{i=1}^{\infty}(\mathrm{diam} (U_i))^s: \bigcup_i U_i\supset F\right\}.
\]
The Hausdorff dimension of $F$ is:
\[
\Haus F=\inf\{s\geq 0:\mathcal{H}^s_\infty(F)=0\}.
\]

\subsection{Box dimensions}
We let $N_r(F)$ be the minimal covering number of a bounded set $F$ in $\mathbb{R}^n$ by cubes of side length $r>0$. The upper/lower box dimension of $F$ is:
\[
\uboxd \text{ resp. }\lboxd(F)=\limsup_{r\to 0} \text{ resp. }\liminf_{r\to 0}\left(-\frac{\log N_r(F)}{\log r}\right).
\]
If the limsup and liminf are equal we call this value the box dimension of $F$.

\subsection{The Assouad dimension}
The \textit{Assouad dimension} of $F$ is
\begin{align*}
\Assouad F = \inf \Bigg\{ s \ge 0 \, \, \colon \, (\exists \, C >0)\, (\forall & R>0)\,  (\forall r \in (0,R))\, (\forall x \in F) \\
&N_r(B(x,R) \cap F) \le C \left( \frac{R}{r}\right)^s \Bigg\}
\end{align*}
where $B(x,R)$ denotes the closed ball of centre $x$ and radius $R$.

In general we have $\Haus(F)\le \lboxd(F)\le \uboxd(F) \le \Assouad(F)$, and all inequalities can be strict. Lower and upper box dimensions and Assouad dimension are invariant under taking closures.

\subsection{Notational conventions}

When counting covering numbers, it is convenient to introduce notions $\approx, \lessapprox, \gtrapprox$ for approximately equal, approximately smaller and approximately larger. The letter $0<\delta\ll 1$ will denote a small scale. Then for two quantities depending on the scale $\delta$, say $f(\delta)$ and  $g(\delta)$ we define the following:
\[
f\lessapprox g\iff \forall \eps>0, \exists C_{\eps}>0 \text{ such that } \forall \delta>0, f(\delta)\leq C_{\eps} \delta^{-\eps} g(\delta).
\]
\[
f\gtrapprox g\iff g\lessapprox f.
\]
\[
f\approx g\iff f\lessapprox g \text{ and } g\lessapprox f.
\]
The constants $C_\eps$ may depend on other variables, so long as these remain independent of $\delta$.

We also use the notation $A\lesssim B$ for $A\le C B$ where $C>0$ is a constant independent of $\delta$ (it may depend on other fixed parameters, so long as they remain independent of $\delta$). Likewise we define $B\lesssim A$ and $A\sim B$.

\section{Lower bound: Planar dipole Kakeya sets}

In this section we prove the lower bound
\[
\lboxd A\geq 4/7
\]
valid for all sets $A\in\DK_2$. Let $\delta>0$ be a small number (<$0.0001$). Let $E_\delta\subset S^1$ be a maximal $\delta$-separated set. In particular, $\#E_\delta\geq 1/(2\delta)$. For each $e\in E_\delta$, there is a pair $(x,y)\in A\times A$ such that $y-x=e$. We collect all such pairs of points for all $e\in E_\delta$. In this way, we obtain a set $A_\delta\subset A$ with cardinality at most $2\#E_\delta$.

Note that $A_\delta$ needs not be $\delta$-separated, so instead of dealing with it directly we study its intersection with $\delta$-squares. We cover $\mathbb{R}^2$ by axes-parallel squares of side length $\delta$ with disjoint interiors, and denote by $\mathcal{C}_\delta$ the collection of all such squares that intersect $A_\delta$. Then $\#C_\delta\sim N_{\delta}(A_\delta)\le N_\delta(A)$, so our task is to estimate
$\#C_\delta$ from below.

Given $C\in \mathcal{C}_\delta$, we write $\dir C$ for the set of directions associated with $C$. More precisely, we define
\be \label{eq:def-dir}
\dir C=\{e\in E_\delta: \exists x\in C, y\in A_\delta,  y-x=e\}.
\ee
By construction,
\[
E_\delta\subset \bigcup_{C\in \mathcal{C}_\delta} \dir C.
\]
Let $\gamma\in (0,1/2)$ be a parameter to be fixed at the end of the proof. We define a set of good directions (depending on $\gamma$):
\[
\good(\gamma) = \big\{ e\in E_\delta: \exists C\in\mathcal{C}_\delta \text{ with } \#\dir C\cap B_{\delta^{1/2}}(e)\geq \delta^{-\gamma} \big\}.
\]
Here $B_{\delta^{1/2}}(e)$ is the $\delta^{1/2}$-ball centred at $e$ in $S^1$.  
We split the argument depending on whether most directions are good or not.

\subsection{Case (i): $\# \good(\gamma) \geq  \# E_\delta/2$}
	
	In this case, we shall use the $L^2$ Kakeya maximal inequality in $\mathbb{R}^2$. Let $\mathcal{K}_\delta$ be the Kakeya maximal operator in $\mathbb{R}^2$; namely, for each function $f\in L^2(\mathbb{R}^2)$ we write
	\[
	s\in S^1\to\mathcal{K}_\delta f(s)=\sup_x \frac{1}{\delta}\int_{T_\delta(s,x)} f(y)dy.
	\]
	Here $T_\delta(s,x)$ is the $1\times \delta$ rectangle (tube) centred at $x\in\mathbb{R}^2$ whose long side has direction $s$. Cordoba's maximal Kakeya inequality (see \cite{Cordoba77}) states that
	\be \label{*}
	\|\mathcal{K}_{\delta}f\|_2\lesssim \sqrt{\log (1/\delta)} \|f\|_2.
	\ee
	
	Now we rescale everything by a factor $\delta^{-1/2}$. Let $e\in\good(\gamma)$. Then, by construction, we can find a $2\delta^{1/2}\times 2\delta$ rectangle $T_e$ whose long side has direction $e^{\perp}$, and such that $T_e$ contains at least $\delta^{-\gamma}/100$ squares in $\mathcal{C}_\delta$. After the rescaling, $T_e$ becomes a $2\times 2\delta^{1/2}$ rectangle and $\mathcal{C}_\delta$ becomes a collection of $\delta^{1/2}$ squares. We write $\tilde{A}_\delta,\tilde{\mathcal{C}}_\delta$ for the rescaled versions of $A_\delta, \mathcal{C}_\delta$. Then there are at least $\#E_\delta/2$ many $\delta$-separated directions, such that for each one of these directions, say $e$, there is a $2\times 2\delta^{1/2}$-rectangle in direction $e$ that contains at least $\delta^{-\gamma}$ many squares in $\tilde{\mathcal{C}}_\delta$. Let $f$ be the characteristic function of the union of squares in $\tilde{\mathcal{C}}_\delta$. Then for $s$ in a $(\delta^{1/2}/100)$-neighbourhood of $e^{\perp},e\in \good(\gamma)$ we have
\[
\mathcal{K}_{\delta^{1/2}} f(s) \gtrsim \frac{1}{\delta^{1/2}} \delta^{-\gamma} (\delta^{1/2})^2 = \delta^{1/2-\gamma},
\]
and therefore, using the assumption $\# \good(\gamma) \geq  \# E_\delta/2\gtrsim \delta^{-1}$,
	\[
	\|\mathcal{K}_{\delta^{1/2}} f\|_2^2\gtrsim \delta^{1- 2\gamma}.
	\]
   On the other hand, since $\#\mathcal{C}_\delta=\#\tilde{\mathcal{C}}_\delta$,
	\[
	\|f\|_2^2\lesssim \delta\#\mathcal{C}_\delta.
	\]
	Thus by the Kakeya maximal inequality $\eqref{*}$ we conclude that
	\be \label{eq:first-case}
	\#\mathcal{C}_\delta\gtrsim \delta^{-2\gamma}/\log(1/\delta)\gtrapprox \delta^{-2\gamma}.
	\ee
	\subsection{Case (ii): $\# \good(\gamma) <  \# E_\delta/2$}
For simplicity, let us define
\[
\bad(\gamma) = \big\{ e\in E_\delta: \forall C\in\mathcal{C}_\delta, \#\dir C\cap B_{\delta^{1/2}}(e)< \delta^{-\gamma} \big\} = E_\delta\setminus \good(\gamma).
\]

Now we construct a graph as follows: the vertex set is $\mathcal{C}_\delta$, and there is an edge between $C_1,C_2\in \mathcal{C}_\delta$ if and only if there are points $x_1\in C_1\cap A_\delta, x_2\in C_2\cap A_\delta$ such that
		\[
		x_2-x_1\in \bad(\gamma).
		\]
		In this case, we write $C_1\sim C_2$. Note that
\[
\#\{ \text{edges in the graph}\} \gtrsim \#\bad(\gamma)\gtrsim \delta^{-1}.
\]
Let $\mathcal{C}'_\delta$ denote the set of vertices of degree at most $100\delta^{-\gamma}$ in the graph. Then the number of edges adjacent to some vertex in $\mathcal{C}'_\delta$ is $\lesssim   \delta^{-\gamma}\#\mathcal{C}_\delta$. Now this implies that there are two possibilities. Let $M$ be a large constant. Either
\be \label{eq:second-case-a}
\#\mathcal{C}_\delta \ge M^{-1}(1/\delta)^{1-\gamma},
\ee
or (provided $M$ is taken large enough) the number of edges adjacent to a vertex in $\mathcal{C}'_\delta$ is at most half the number of total edges in the graph, which we recall is $\gtrsim 1/\delta$. Suppose from now on we are in the latter case, that is
\be \label{eq:edges-in-graph}
\sum_{C\in\mathcal{C}_\delta\setminus \mathcal{C}'_\delta} \#\{C'\sim C\} \gtrsim \delta^{-1}.
\ee
	
Given two squares $C,C'$, we let $d(C,C')$ denote the distance between their centres. Note that if $C_1,C_2$ satisfy $d(C_1,C_2) \ge \delta^{1/2}$, then
\be  \label{eq:double-degree}
\#\{ C\in \mathcal{C}_\delta: C\sim C_1, C\sim C_2\} \lesssim (1/\delta)^{\gamma}.
\ee
Indeed, the centre of such $C$ must be contained in the intersection of two annuli $A_i$ centred at the centre of $C_i$ with unit radius and width $10\delta$. The intersection of these two annuli is contained in the union of at most two $100\delta^{1/2}$-arcs of $A_i, i\in\{1,2\}$. To see this, consider two unit circles $L_1,L_2$ whose centres have distance at least $\delta^{1/2}$ and at most $(2-\delta^{1/2})$; note that the two circles intersect at two points. At each of those two intersections, the unit tangent vectors of $L_1,L_2$ are $\gtrsim\delta^{1/2}$-separated.  Thus, the points on $L_1$ at distance $\leq \delta$ from $L_2$ can be covered by a $(2\delta)\times(100\delta^{1/2})$-rectangle. If the centres of $L_1,L_2$ are at distance at least $(2-\delta^{1/2})$  from each other, then the two circles are almost outer tangent. In this case, it is possible to see that the intersection of the $\delta$-neighbourhoods of $L_1,L_2$ can be covered by one $(2\delta)\times(100\delta^{1/2})$-rectangle. Recalling the definitions of $\bad(\gamma)$ and $\sim$, we see that \eqref{eq:double-degree} holds.
		
Now, given $C\in \mathcal{C}_\delta$, let
\[
P(C) =\{ (C_1,C_2)\in \mathcal{C}_\delta^2: C\sim C_1, C\sim C_2, d(C_1,C_2) \ge \delta^{1/2}\}.
\]
We claim that if $C\in\mathcal{C}_\delta\setminus \mathcal{C}'_\delta$, then
\be \label{eq:deg-squared}
\#P(C)\gtrsim \deg(C)^2 = \#\{C'\sim C\}^2.
\ee
Indeed, it is enough to show that for any $C_1\sim C$, there are $\gtrsim\deg(C)$ neighbours $C_2\sim C$ such that $d(C_1,C_2)\ge \delta^{1/2}$. This holds since the $C'$ such that $C'\sim C$, $d(C_1,C')<\delta^{1/2}$ are contained in the $4\delta$-neighborhood of a unit circle arc of length $2\sqrt{\delta}$ centred at $x_0\in C$ and therefore, by the definition of $\bad(\gamma)$ and the graph, they number less than $40 \delta^{-\gamma}<\deg(C)/2$, since we assumed that $C\in\mathcal{C}_\delta\setminus \mathcal{C}'_\delta$.

Putting these facts together and using Cauchy-Schwarz, we conclude that
\begin{align*}
\delta^{-\gamma} (\#\mathcal{C}_\delta)^2  &\overset{\eqref{eq:double-degree}}{\gtrsim}  \sum_{C\in\mathcal{C}_\delta\setminus\mathcal{C}'_\delta} \#P(C) \overset{\eqref{eq:deg-squared}}{\gtrsim}  \sum_{C\in\mathcal{C}_\delta\setminus\mathcal{C}'_\delta} \#\{C'\sim C\}^2 \\
&\overset{\text{C-S}}{\gtrsim} \frac{1}{\#\mathcal{C}_\delta} \left(\sum_{C\in\mathcal{C}_\delta\setminus\mathcal{C}'_\delta} \#\{C'\sim C\}\right)^2\overset{\eqref{eq:edges-in-graph}}{\gtrsim}  \frac{\delta^{-2}}{\#\mathcal{C}_\delta},
\end{align*}
and hence
\be \label{eq:second-case-b}
\#\mathcal{C}_\delta \gtrsim (1/\delta)^{(2-\gamma)/3}.
\ee

We have finished our discussion. To summarise, combining \eqref{eq:first-case}, \eqref{eq:second-case-a} and \eqref{eq:second-case-b} we get that
\[
\#\mathcal{C}_\delta\gtrapprox (1/\delta)^{\min(2\gamma,1-\gamma,(2-\gamma)/3)}.
\]
At this point we observe that the optimal choice is $\gamma=2/7$ and we conclude that
\[
N_\delta(A) \gtrsim \#\mathcal{C}_\delta\gtrapprox (1/\delta)^{4/7},
\]
leading to
\[
\lboxd A\geq 4/7.
\]

We remark that this proof appears to break down in dimension $n\ge 3$. The issue is that the intersection of $n$ (or even more) $\delta$-neighborhoods of spheres with $\delta^{1/2}$-separated centres needs not be contained in a small number of $\delta^{1/2}$-balls. Hence there seems to be no useful analogue of \eqref{eq:double-degree}. Thus it remains an open problem to give non-trivial lower bounds for the box dimension of dipole Kakeya sets in dimension $n\ge 3$.

\subsection{Assouad dimension}

We shall now prove that $\Assouad A\geq 2/3$. The argument is somewhat similar to that of case (ii) above. Let $S'\subset S^1$ be an arc of angle $\pi/10$. We consider the subset $A'$ of $A$ formed by dipoles pointing in directions in $S'$. Let $E_\delta$ be a maximal $\delta$-separated subset of $S'.$ We define $\mathcal{C}_\delta$ as above. Given $C_1,C_2\in \mathcal{C}_\delta$, we write $C_1\sim C_2$ if there are points $x_1\in C_1, x_2\in C_2$ with $x_2-x_1\in E_\delta$. Observe that if $C_1\sim C_2$ and $C_1\sim C_3,$ then we have $d(C_2,C_3)<1.9$.

Now consider the incidence set
\[
\mathcal{I}=\{(C_1,C_2,C_3)\in \mathcal{C}^3_\delta: C_1\sim C_2,C_1\sim C_3 \}.
\]
For each $C_1\in\mathcal{C}_\delta,$ let $n_{C_1}$ be the number of $C_2\in\mathcal{C}_\delta$ with $C_1\sim C_2.$ Then we see that
\begin{align}\label{eq:cs}
\#\mathcal{I}\gtrsim \sum_{C_1\in\mathcal{C}_\delta} n^2_{C_1}\geq \frac{1}{\#\mathcal{C}_\delta}\left(\sum_{C_1\in\mathcal{C}_\delta} n_{C_1}\right)^2\gtrsim \frac{1}{\#\mathcal{C}_\delta}\delta^{-2}.
\end{align}
We have used Cauchy-Schwarz for the second inequality above.  Now fix
\[
t > \Assouad A.
\]
By the definition of Assouad dimension, there is a constant $M_t$ (independent of $\delta$) such that for each $C\in \mathcal{C}_\delta$ and integer $j\geq 0$,
\[
\#\big\{ C_1\in\mathcal{C}_\delta: d(C,C_1)\leq 2^{j+1}\delta  \big\} \le M_t (2^{j}\delta/\delta)^{t}.
\]
This implies that the number of pairs $C_1,C_2\in\mathcal{C}_\delta$ with $d(C_1,C_2)\leq 2^{j+1}\delta$ is $\lesssim M_t2^{tj}\#\mathcal{C}_\delta.$ Let $(C_1,C_2)$ be such that $d(C_1,C_2)\in [2^j\delta, 2^{j+1}\delta)$. We want to estimate the number of $C_0\in\mathcal{C}_\delta$ such that
\[
(C_0,C_1,C_2)\in \mathcal{I}.
\]
We have considered this problem above when $C_1,C_2$ are sufficiently separated. By our above discussion, we may assume $2^j\delta\in (0,1.9)$. The point of this condition is to avoid 'outer-tangencies' between two unit circles. More precisely, we consider the intersection between two annuli with radius $1$ and width $\delta$. The shape of this intersection depends on the distance $d$ between the centres of these two annuli. If $d<1.9$ then the intersection is contained in a $10\delta$ neighbourhood of two $\lesssim \delta/d$ arcs of one of the annuli. However, when $d$ is near $2$, then the two annuli are 'outer-tangent' and their intersection is roughly a $\delta\times \delta^{1/2}$ rectangle which is too large for our purposes.

Now the number of $C_0$ with $(C_0,C_1,C_2)\in \mathcal{I}$ can be bounded from above by
\[
M_t\left(\frac{\delta/2^j\delta}{\delta}\right)^{t}
\]
for a (different) constant $M_t>0$. Then we see that
\be \label{+++}
\#\mathcal{I}\leq M_t\sum_{j}  M_t2^{tj}\#\mathcal{C}_\delta \left(\frac{\delta/2^j\delta}{\delta}\right)^{t}\lesssim M^2_t (1/\delta)^{t} \#\mathcal{C}_\delta \log(1/\delta).
\ee
Combining \eqref{eq:cs}, \eqref{+++} we deduce
\[
\delta^{-2}\lesssim (\#\mathcal{C}_\delta)^2 \delta^{-t}\log (1/\delta)\lesssim \delta^{-3t}\log (1/\delta).
\]
Thus $t\ge 2/3$, and since $t>\Assouad A$ was arbitrary this concludes the proof.

\section{Constructions}\label{RE}
\subsection{Construction for lower box dimension} \label{subsec:construction-lower-box}
\begin{figure}[H]
	\includegraphics[width=0.8\textwidth]{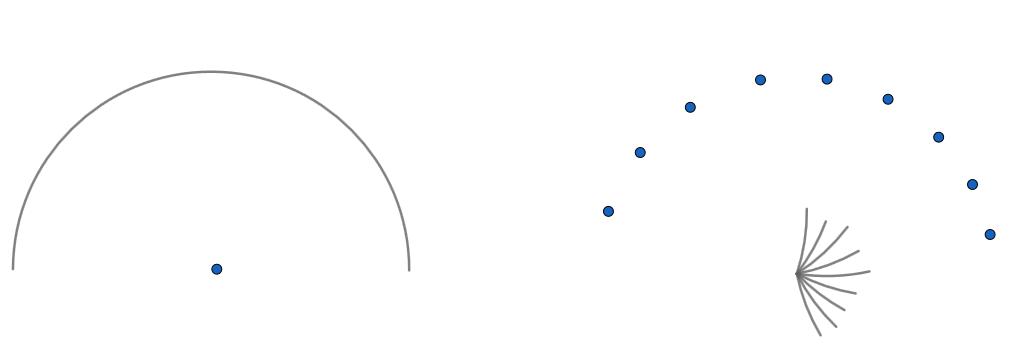}
	\caption{Operation on unit circular arcs}
	\label{arc2}
\end{figure}

We construct a dipole Kakeya set $P\in \DK_2$ with lower box dimension at most $2/3$. The construction is based on the following operation acting on arcs with unit radius. This operation is illustrated in Figure \ref{arc2}.

Given an arc $E$ with unit radius centred at $e\in\mathbb{R}^2$ and a positive number $\delta>0$, we define two sets $\mathcal{P}(E,\delta)$ and $\mathcal{A}(E,\delta)$ (here $\mathcal{P}$ stands for ``points'' and $\mathcal{A}$ for ``arcs'') in the following way:
First we choose a partition of $E$ into sub-arcs with length between $\delta$ and $2\delta$. Such a partition exists when the arc length of $E$ is larger than $\delta$. If this is not the case, we use the trivial partition which is $E$ itself. In general, multiple choices of such partitions are possible, we only need to choose one. Suppose this partition we chose is determined by the following sequence of points on $E$ in the clockwise order:
\[
e_1,e_2,e_3\dots, e_N.
\]
Then we define $\mathcal{P}(E,\delta)=\{e_i\}_{i=1,\dots,N}$. Now for any two neighbouring points $e_i,e_{i+1}$, we consider the arc $e_ie_{i+1}$ between them. We want to 'transfer' this arc to the point $e$ in such a way that the unit directions determined by the set are preserved. To do this, we draw an arc $ee'_{i}$ with centre $e_{i}$ starting from $e$ in counterclockwise direction such that the arc length of $ee'_i$ is equal to that of $e_{i}e_{i+1}$. We let $\mathcal{A}(E,\delta)$ be the collection of all such arcs
\[
\mathcal{A}(E,\delta)=\{ee'_{i}\}_{i=1,\dots,N-1}.
\]
We choose a sequence of positive numbers which decays very fast, for example,
\[
\delta_k=\frac{1}{2^{2^{(k^2)}}}.
\]
First take $E_0$ to be the unit circle $S^1$ and $P_0=\{0\}$, and apply the above procedure to obtain $\mathcal{P}(E_0,\delta_1)$, $\mathcal{A}(E_0,\delta_1)$. We denote
\[
P_1=\mathcal{P}(E_0,\delta_1),\quad A_1=\mathcal{A}(E_0,\delta_1).
\]
In general after the $k$-th step we have the sets $P_k, A_k$ of points and arcs respectively. Then, for each element $E\in A_k$ we obtain $\mathcal{P}(E,\delta_{k+1})$, $\mathcal{A}(E,\delta_{k+1})$. Then denote
\begin{align*}
P_{k+1}&=\bigcup_{E\in A_k}\mathcal{P}(E,\delta_{k+1}),\\
A_{k+1}&=\bigcup_{E\in A_k}\mathcal{A}(E,\delta_{k+1}).
\end{align*}
The idea of this construction is that if the union of the points and the arcs in $P_k, A_k$ contains a unit distance in every direction, so does the union of the points and the arcs in $P_{k+1}, A_{k+1}$, while the latter set is much smaller because many arcs have moved closer together.

Then we define $P=\bigcup_{k\in\mathbb{N}} P_k$. For each integer $k$, the set $\bigcup_{i\leq k} P_i$ contains two endpoints of unit segments whose directions form a $(2\delta_{k})$-dense subset of $S^1$. Then we take the closure $\overline{P}$ and finish our construction. Since $\overline{P}-\overline{P}$ contains a dense subset of $S^1$ by our previous observation, it contains all of $S^1$. Our task is then to get an upper bound for the lower box dimension of $\overline{P}$. To begin, as the lower box dimension is stable under taking closures, it is enough to show that $\lboxd P\le 2/3$. To do this, we consider the scales $r_k=\delta^{3/2}_k$ and covering numbers $N_{r_k}(P)$,  for large enough $k\in\mathbb{N}$. In order to estimate $N_{r_k} (P)$ it is enough to compute $N_{r_k}\left(\bigcup_{i\leq k+1} P_i\right)$, because $P$ is contained in an $\eps$-neighbourhood of $\bigcup_{i\leq k+1} P_i$, where $\eps$ can be bounded by
\[
\eps\leq\sum_{i\geq k+1} \delta_{i}\leq 2\delta_{k+1}\leq (\delta_k)^{3/2} = r_k,
\]
using the fast decrease of $\delta_k$. Hence
\be  \label{@}
N_{r_k}\left(\bigcup_{i\leq k+1} P_i\right)\sim N_{r_k}\left(P\right).
\ee
We estimate the box covering number of $\bigcup_{i\leq k+1} P_i$ with respect to the scale $r_k$. For $\bigcup_{i\leq k} P_i$, we see that the cardinality of this set is at most $\lesssim \sum_{i\leq k} \delta^{-1}_{i}$.

Now we make the key observation that $P_{k+1}$ is contained in a $(2\delta_k)$-neighbourhood of $P_{k-1}$ by construction, and hence
\[
N_{r_k}(P_{k+1}) \lesssim (\delta_k/r_k)^2 \#P_{k-1} \lesssim \delta_k^{-1}\delta_{k-1}^{-1}.
\]
Therefore, using again the fast decay of $\delta_{k}=\frac{1}{2^{2^{(k^2)}}}$,
\[
N_{r_k}\left(\bigcup_{i\leq k+1} P_i\right)\lesssim \delta_k^{-1}\delta_{k-1}^{-1} +\sum_{i\leq k} \delta^{-1}_{i} \lesssim \delta^{-1}_{k-1} \delta^{-1}_k.
\]
Note that $\log(\delta_{k-1})/\log(\delta_k)\to 0$ by the fast decay of $\delta_k$.  Recalling \eqref{@}, we conclude that
\[
\liminf_{k\to\infty}\frac{\log N_{r_k}(P)}{\log(1/r_k)}\leq \liminf_{k\to\infty} \frac{\log(\delta_{k-1}\delta_{k})}{(3/2)\log \delta_k}=2/3.
\]
We have therefore verified that
\[
\lboxd \overline{P}=\lboxd P\leq \frac{2}{3}.
\]

We can also perform similar constructions in $\mathbb{R}^n$ for $n\geq 3$. The idea is exactly the same as in $\mathbb{R}^2$; we only sketch the construction. We first need to modify the definition of $\mathcal{P}(E,\delta)$ and $\mathcal{A}(E,\delta)$. Here $E$ is an open ball in a unit sphere $S$, say with center $e$. We now let $\mathcal{P}(E,\delta)$ be a maximal $\delta$-separated subset of $E$. For each $p\in \mathcal{P}(E,\delta)$, let $E_p$ be an isometric copy of $E'_p=B(p,2\delta)\cap S$ such that $e\in E_p$, and the set of directions spanned by $p$ and $E_p$ is the set of directions spanned by $e$ and $E'_p$. We let $\mathcal{A}(E,\delta)$ be the union of all $E_p$, $p\in\mathcal{P}(E,\delta)$. Then we can perform the construction of the set $P\subset\mathbb{R}^n$ with a sequence of $\delta_i\to 0$ and estimate the covering numbers at the scales $\delta^{(2n-1)/n}_i$. Arguing as in the planar case, we are led to the bound
\[
\lboxd \overline{P}\leq \frac{n(n-1)}{2n-1}.
\]

\subsection{The Hausdorff dimension of $\overline{P}\subset\mathbb{R}^2$}
\label{subsec:Hausdorff}

The computation of $\Haus \overline{P}$ is not very different from that of $\lboxd \overline{P}$. The key is again that $P_{k+1}$ is contained in a $(2\delta_k)$-neighbourhood of $P_{k-1}$, and together with the fast decay of $\delta_k$ this implies that
\[
\overline{P} \subset \left(\bigcup_{j=0}^{k-2} P_j\right) \cup P_{k-1}^{(3\delta_k)} \cup P_k^{(3\delta_{k+1})},
\]
where $F^{(\rho)}=\{ x: \dist(x,F)<\rho\}$. Therefore we can cover $\overline{P}$ by $\sum_{j=0}^{k-1} \#P_j$ balls of radius $3\delta_k$ and $\#P_k$ balls of radius $3\delta_{k+1}$. Hence, for any $s>0$ we can bound
\[
\mathcal{H}_\infty^s(\overline{P}) \lesssim \big(\sum_{i\le k-1} \delta_i^{-1}\big) \delta_k^s + \delta_k^{-1} \delta_{k+1}^s.
\]
Using the fast decay of $\delta_i$ we see that the right-hand side tends to $0$ as $k\to\infty$, and hence $\mathcal{H}_\infty^s(\overline{P})=0$ and $\Haus \overline{P}=0$. The same argument works in arbitrary dimension.

\subsection{A construction with upper box dimension equal to $3/4$}
We can also construct a planar dipole Kakeya set with upper box dimension equal to $3/4$. The idea is not very different from the construction above. We provide some details.

\begin{figure}[H]
	\includegraphics[width=\textwidth]{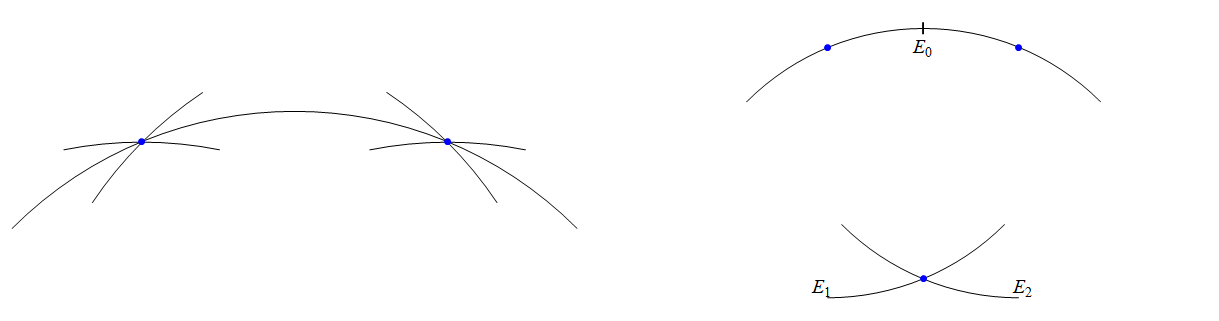}
	\caption{The operation on arcs (left) and the first step of the construction (right)}
	\label{arc3}
\end{figure}

Let $E$ be an arc of a unit circle. We will define collections of arcs $\mathcal{A}(E)$ and of points $\mathcal{P}(E)$ depending on $E$. Let $\theta\in (0,\pi)$ be the angle of $E$. Let $x_1,x_2$ be the two points, in addition to the midpoint of $E$, needed to cut $E$ into four equal pieces. Let $E_i$ be the sub-arc of $E$ of angle $\theta/4$ having $x_i$ as its midpoint.  Let $E'_i,E''_i$ be two rotated copies of $E_i$ around $x_i$ by angles $\pm \theta/8$ respectively. See Figure \ref{arc3}. Note that the centres of these arcs are obtained by rotating the centre of $E$ around $x_i$ by $\pm\theta/8$. (Here and below, by the centre of an arc we mean the centre of the circle that contains it, not to be confused with the midpoint of the arc.) We define $\mathcal{A}(E)=\{E'_1,E''_1,E'_2,E''_2\}$, and $\mathcal{P}(E)=\{x_1,x_2\}$. Starting with $E$, we inductively perform these operations on all the sub-arcs obtained in the process, and define $F(E)$ to be the closure of the countable set of all the points in $\mathcal{P}(E')$ obtained along the way (over all arcs $E'$ obtained at any step of the process). Our final set $F$ is defined as the union of $F(E_0), F(E_1), F(E_2)$  where $E_0$ is a $\pi/2$-arc and $E_1,E_2$ are $\pi/4$ arcs as depicted in Figure \ref{arc3}. Namely, $E_1$ and $E_2$ intersect at the centre of the arc $E_0$, and the centres of $E_1$ and $E_2$ are the points $x_1$, $x_2$ associated to the arc $E_0$. Hence this configuration consists of two pairs of ``opposite'' $\pi/4$-arcs.
By the choice of the initial arcs in ``opposite'' configuration, at each stage in the process the unions of points and arcs contains every unit direction. Hence the same argument from \S\ref{subsec:construction-lower-box} shows that $S^1\subset  (F \cup F')-(F\cup F')$, where $F'$ is a copy of $F$ rotated by $90$ degrees.

Since upper box dimension is finitely stable,
\[
\uboxd F=\max\{\uboxd F(E_0),\uboxd F(E_1),\uboxd F(E_2)\}.
\]
We now show that $\uboxd F(E_0)=3/4$; the same argument applies to $E_1$ and $E_2$ as well.

Let $k$ be a large integer. We start with $E_0$ and perform the above operation until all the smaller arcs are of angle $4^{-k}\pi/2.$ As a result, the collection $P_k$ of pairs of points at this stage has cardinality $\lesssim 4^k$. We also have $4^k$ arcs of angle $4^{-k}\pi/2.$ We want to count how many balls of radius $4^{-k}$ are needed to cover all the points in $P_k$. If we then enlarge all the balls tenfold, the union will also cover all the arcs.

Denote this counting number by $N_k.$ It equals $N(F(E_0),10\cdot 4^{-k})$ up to fixed multiplicative constants, and so tells us everything about the upper box dimension of $F(E_0)$. A very naive estimate for $N_k$ is simply $100\cdot 4^k.$ However, the points in $P_k$ are not $4^{-k}$-separated. Indeed, $P_{n+1}$ is constructed by adding four points around each element of $P_n$. Those four points are the pairs of endpoints of two arcs of angle $4^{-n-1}\pi/4$, and those arcs have tangential directions at their midpoint separated by an angle $4^{-n-1}\pi/2$. Thus, these four points can be grouped into two pairs, say, $x_1,x_2$ and $y_1,y_2$, such that $d(x_1,x_2)=d(y_1,y_2) \sim 4^{-2n}$ and $d(x_i,y_j)\sim 4^{-n}$. Recall we want to estimate the $4^{-k}$-covering number. The idea is now that for $n=\lfloor k/2+100\rfloor$ to $n=k$, when seen at the resolution $4^{-k}$, each point in $P_n$ splits with multiplicity $2$ rather than $4$. To be more precise, what holds is that if $x\in P_{\lfloor k/2+100\rfloor}$, then all points that split from $x$ in all construction steps from $n=\lfloor k/2+100\rfloor$ to $n=k$ can be covered by $2^{1+k-\lfloor k/2+100\rfloor}$ balls of radius $4^{-k}$. Hence
\[
N_k \le 4^{1+\lfloor k/2+100\rfloor} 2^{1+k-\lfloor k/2+100\rfloor}\lesssim 4^{k/2}2^{k/2} = 4^{3k/4}.
\]
In fact, the same argument shows that $N_k\approx 4^{3k/4}$. As explained above, this leads to $\uboxd F(E_0)=3/4$, which is the result we claimed. (In fact, the argument shows that $\lboxd F(E_0)=3/4$ as well.)

\section*{Acknowledgements} PS thanks Alan Chang for suggesting the problem, and Marianna Cs\"{o}rnyei and Tam\'{a}s Keleti for useful discussions during the early phase of the project. HY thanks Korn\'{e}lia H\'{e}ra for useful discussions. HY thanks the Corpus Christi College, Cambridge for financial support. Both authors thank an anonymous referee for pointing out an inaccuracy in an earlier version of the manuscript, suggesting a simplification for the proof of the lower bound for Assouad dimension, and other helpful remarks.

\begin{dajauthors}
\begin{authorinfo}[pablo]
Pablo Shmerkin\\
Department of Mathematics and Statistics\\
Torcuato Di Tella University, and CONICET\\
Av. Figueroa Alcorta 7350 (1428)\\
Buenos Aires, Argentina\\
Current address: Department of Mathematics\\
The University of British Columbia\\ 
1984 Mathematics Road, Vancouver BC V6T 1Z2, Canada \\
pshmerkin\imageat{}utdt\imagedot{}edu\\
\url{https:pabloshmerkin.org}
\end{authorinfo}
\begin{authorinfo}[han]
 Han Yu \\
  Department of Pure Mathematics and Mathematical Statistics\\
  University of Cambridge\\
  CB3 0WB, UK\\ 
  hy351\imageat{}cam\imagedot{}ac\imagedot{}uk
\end{authorinfo}
\end{dajauthors}

\end{document}